\documentclass[11pt]{article}
\usepackage{amsfonts}
\usepackage{eufrak}
\usepackage{mathrsfs}
\usepackage{amssymb}
\usepackage{mathtools}
\usepackage{commath}
\usepackage{dsfont}
\usepackage{amsthm,amsmath}
\begin{document}
\title{\bf The exceptional set for Diophantine approximation \\
with mixed powers of prime variables
\thanks{Project supported by the National Natural Science Foundation of China (grant no.11771333) and by the Fundamental Research Funds for the Central Universities (grant no. JUSRP122031).}}
\author{{Yuhui  Liu}\\ {School of Science,  Jiangnan University}\\{Wuxi 214122, Jiangsu, China}
\\{Email: tjliuyuhui@outlook.com}}
\date{}
\maketitle
{\bf Abstract} Let $\lambda_1, \lambda_2, \lambda_3, \lambda_4$ be non-zero real numbers, not all negative, with $\lambda_1/\lambda_2$ irrational and algebraic. Suppose that $\mathcal{V}$ is a well-spaced sequence and $\delta >0$.
In this paper, it is proved that  for any $\varepsilon >0$, the number of $v \in \mathcal{V}$ with $v \leqslant N$ for which
\begin{align*}
|\lambda_1 p_1^2 + \lambda_2 p_2^3+ \lambda_3 p_3^4+ \lambda_4 p_4^5 - v| < v^{-\delta}
\end{align*}
has no solution in prime variables $p_1,p_2,p_3,p_4$ does not exceed $O\big(N^{\frac{359}{378} + 2\delta +\varepsilon}\big)$. This result constitutes an improvement upon that of Q. W. Mu and Z. P. Gao [12].

{\bf 2010 Mathematics Subject Classification}: 11D75, 11P32.\

{\bf Key Words}: Diophantine inequality; prime; Davenport-Heilbronn method; Exceptional set.

\section{Introduction}
\setcounter{equation}{0}
\hspace{5.8mm} Let $\lambda_1, \lambda_2, \lambda_3, \lambda_4$ are nonzero real numbers, not all negative. Suppose that $\lambda_1/\lambda_2$ is irrational. In 2015, Yang ang Li [17] considered the inequality
 \begin{align}
\left|\lambda_1 x_1^2 + \lambda_2 x_2^3+ \lambda_3 x_3^4+ \lambda_4 x_4^5 - p - \frac12\right| < \frac12,
\end{align}
and proved that the inequality (1.1) has infinite solutions with prime $p$ and $x_1, x_2, x_3, x_4 \in \mathbb{Z}$.

It is of some interest to investigate the exceptional set for the analogous inequality
\begin{align}
|\lambda_1 p_1^2 + \lambda_2 p_2^3+ \lambda_3 p_3^4+ \lambda_4 p_4^5 - v| < v^{-\delta}.
\end{align}
An increasing sequence $v_1 < v_2 < \cdot\cdot\cdot$ of positive real numbers is called a well-spaced sequence if there exist positive constants $C > c > 0$ such that
\begin{align*}
0 < c < v_{i+1} - v_i < C, \,\,\mbox{for}\,\, i = 1, 2, ....
\end{align*}
In 2018, Motivated by  Harman [8] and Languasco and Zaccagnini [11], Ge and Zhao [7] established the following Theorem:

{\bf Theorem 1.0} {\it \,\,Suppose that $\lambda_1, \lambda_2, \lambda_3, \lambda_4$ are nonzero real numbers, not all negative, $\lambda_1/\lambda_2$ is irrational and algebraic, $\mathcal{V}$ is a well-spaced sequence, $\delta >0$. Let $E(\mathcal{V}, N, \delta)$ denote the number of $v \in \mathcal{V}$ with $v \leqslant N$ for which
 \begin{align}
|\lambda_1 p_1^2 + \lambda_2 p_2^3+ \lambda_3 p_3^4+ \lambda_4 p_4^5 - v| < v^{-\delta}
\end{align}
has no solution in primes $p_1, p_2, p_3, p_4$. Then for any $\varepsilon >0$},
\begin{align*}
E(\mathcal{V}, N, \delta) \ll N^{1-\frac{1}{72}+2\delta+\varepsilon}.
\end{align*}
Very recently, by using some ingredients due to Br\"{u}dern [2] and Chow [3], the result was improved to
\begin{equation*}
E(\mathcal{V}, N, \delta) \ll X^{1-\frac{13}{360}+2\delta+\varepsilon}
\end{equation*}
by Mu and Gao [12]. In this paper, we obtain a better upper bound of $E(\mathcal{V}, N, \delta)$ by giving the following theorems.\\

{\bf Theorem 1.1} {\it \,\,Subject to the conditions in Theorem 1.0, we have
\begin{align*}
E(\mathcal{V}, N, \delta) \ll N^{1-\frac{19}{378}+2\delta+\varepsilon}
\end{align*}
for any $\varepsilon > 0$.}\\

{\bf Theorem 1.2} {\it \,\,Let $\lambda_1, \lambda_2, \lambda_3, \lambda_4$ be nonzero real numbers, not all negative. Suppose that $\lambda_1/\lambda_2$ is irrational. Let $\mathcal{V}$ be a well-spaced sequence and $\delta >0$. Then there exists a sequence $N_j \rightarrow \infty$ such that for any $\varepsilon >0$,
\begin{align}
E(\mathcal{V}, N_j, \delta) \ll N_j^{1-\frac{19}{378}+2\delta+\varepsilon}.
\end{align}
Moreover, if the convergent denominators $q_j$ for $\lambda_1/\lambda_2$ satisfy that for some $\omega \in [0,1)$,
\begin{align}
q^{1-\omega}_{j+1} \ll q_j.
\end{align}
Then for all $N\geqslant 1$ and any $\varepsilon >0$,
\begin{align}
E(\mathcal{V}, N, \delta) \ll N^{1-19\chi+2\delta+\varepsilon}
\end{align}
with}
\begin{align}
\chi = \min\left(\frac{1-\omega}{94-75\omega},\frac{1}{378}\right).
\end{align}

{\bf Remark}  In the case of $\lambda_1/\lambda_2$ algebraic, Roth's Theorem (see Theorem 2A on page 116 in Schmidt [13]) implies that we can take $\omega=\varepsilon$ in (1.7), and so $\chi = \frac{1}{378}$. Thus Theorem 1.1 follows immediately from Theorem 1.2.
\section{Notation and Some Preliminary Lemmas}
\setcounter{equation}{0}
\hspace{5.8mm}For the proof of the Theorems, in this section we introduce the necessary notation and Lemmas.

Throughout this paper, let $\eta < 10^{-10}$ be a fixed positive constant, and $\varepsilon\in(0, 10^{-10})$ be an arbitrarily small positive constant, not necessarily the same in different formulae. The letter $p$, with or without subscripts, is reserved for a prime number. By $A \sim B$ we mean that $B < A \leqslant 2B$. We use $e(\alpha)$ to denote $e^{2\pi i\alpha}$. $\lambda_1, \lambda_2,\lambda_3,\lambda_4$ are nonzero real numbers , not all of the same sign and $\lambda_1/\lambda_2$ is irrational. Implicit constants in the $O, \ll$ and $\gg$ notations usually depend at most on $\lambda_1, \lambda_2,\lambda_3,\lambda_4, \varepsilon$ and $\eta$. $d(q)$ denotes the divisior function. Since $\lambda_1/\lambda_2$ is irrational, we let $q$ be any denominator of a convergent to $\lambda_1/\lambda_2$, and let $X$ run through the sequence $X=q^{\frac73}$. Let $\tau\in(0, 1)$, and we shall eventually take $\tau = X^{-\delta}$. We define
\begin{eqnarray*}
&&L=\log X,\,\,\,I_j = \bigg[(\eta X)^{\frac {1}{j}}, X^{\frac {1}{j}}\bigg] ,\,\,\,\,K(\alpha) =
\begin{cases}
\left(\frac{\sin(\pi \tau \alpha)}{\pi \alpha}\right)^2,  &\mbox{if}\,\,\,\, \alpha \neq 0;\\
{\tau}^2,  &\mbox{if}\,\,\,\, \alpha = 0.
\end{cases}
\end{eqnarray*}
It is well known that
\begin{align}
K(\alpha) &\ll \min (\tau^2, |\alpha|^{-2}),\\
A(x) &= \int_{-\infty}^{\infty} e(\alpha x)K(\alpha)d\alpha = \max(0, \tau - |x|).
\end{align}
The identity (2.2) is a corollary of Davenport and Heilbronn [4, Lemma 4].

Moreover, we borrow the function $\rho(m)$ defined in (5.2) of Harman and Kumchev [9] (see also Harman [8, Section 8]) and write
\begin{align*}
\psi(m,z) &=
\begin{cases}
1,  &\mbox{if}\,\,\,\,p|m \Rightarrow p\geqslant z;\\
0,  &\mbox{otherwise},
\end{cases}\,\,\,
z(p) =
\begin{cases}
X^{\frac{5}{28}}p^{-\frac12},  &\mbox{if}\,\,\,\,p < X^{\frac17};\\
p,  &\mbox{if}\,\,\,\,X^{\frac17} \leqslant p \leqslant X^{\frac{3}{14}};\\
X^{\frac{5}{14}}p^{-1},  &\mbox{if}\,\,\,\,p > X^{\frac{3}{14}},\\
\end{cases}\\
\rho (m) &= \psi(m, X^{\frac{5}{42}}) - \sum_{X^{\frac{5}{42}} \leqslant p < X^{\frac{1}{4}}}\psi\left(\frac{m}{p}, z(p)\right).
\end{align*}

Let $I^{'}$ be any subinterval of $I_{2}$, by Section 8 in Harman [8], we have
\begin{equation}
\sum_{m \in I^{'}}\rho (m) = \kappa|I^{'}|L^{-1} + O(X^{\frac12}L^{-2}),
\end{equation}
where $\kappa > 0$ is an absolute constant. Furthermore, we define
\begin{eqnarray*}
&&S_2(\alpha)=\sum_{m \in I_2}\rho(m) e(\alpha m^{2}),\,\,\,\,\,\,\,\,\,\,\,\,\,\,T_2(\alpha)=\int_{I_2}e(\alpha t^{2})d\alpha,\\
&&S_k(\alpha)=\sum_{p \in I_k}(\log p) e(\alpha p^{k} )\,\, (k \neq 2),\,\,\,\,\,\,T_k(\alpha)=\int_{I_k}e(\alpha t^{k})d\alpha \,\,(k \neq 2).
\end{eqnarray*}
For any measurable subset $\mathfrak{X}$ of $\mathbb{R}$, write
\begin{eqnarray}
I(\tau, v, \mathfrak{X}) = \int_{\mathfrak{X}}\prod_{j=1}^{4}S_{j+1}(\lambda_{j}\alpha)e(-\alpha v)K(\alpha)d\alpha.
\end{eqnarray}
By (2.2), we have
\begin{align}
I(\tau, v, \mathbb{R}) &= \sum\limits_{{m_{1} \in I_{2}, p_{2} \in I_{3},\atop{p_{3} \in I_{4},p_{4} \in I_{5}}}}\rho(m_1)\prod_{j=2}^{4}(\log p_j) A(\lambda_1 m_1^2 + \lambda_2 p_2^3 + \lambda_3 p_3^4 + \lambda_4 p_4^5 - v)\notag\\
&= \sum\limits_{{m_{1} \in I_{2}, p_{2} \in I_{3},\atop{p_{3} \in I_{4},p_{4} \in I_{5}}}}\rho(m_1)\prod_{j=2}^{4}(\log p_j)\max(0, \tau - |\lambda_1 m_1^2 + \lambda_2 p_2^3 + \lambda_3 p_3^4 + \lambda_4 p_4^5 - v|)
\end{align}
and
\begin{align}
I(\tau, v, \mathbb{R})\ll \tau L^3 \mathcal{N}(X),
\end{align}
where $\mathcal{N}(X)$ denotes the number of solutions to the inequality
\begin{align*}
|\lambda_1 p_1^2 + \lambda_2 p_2^3 + \lambda_3 p_3^4 + \lambda_4 p_4^5 - v| < \tau
\end{align*}
with $p_j \in I_{j+1}$ for $1 \leqslant j \leqslant 4$.

We now divide the real line into the major arc $\mathfrak{M}$, the minor arc ${\mathfrak{m}}$ and the trivial arc ${\mathfrak{t}}$. We define
\begin{align*}
&\mathfrak{M}=\{\alpha :|\alpha|\leqslant P/X\},\,\,\,\mathfrak{m}=\{\alpha : P/X <|\alpha|\leqslant R \},\,\,\,\mathfrak{t}= \{\alpha :|\alpha|> R\},
\end{align*}
where $P=X^{\frac{1}{6}}$ and $R= \tau^{-2}X^{\frac{43}{120} + 2\varepsilon}$. Thus we have
\begin{align}
I(\tau, v, \mathbb{R})= I(\tau, v, \mathfrak{M}) + I(\tau, v, \mathfrak{m}) + I(\tau, v,\mathfrak{t}).
\end{align}

Suppose that $N$ is some large positive quantity which we will choose later. Let
\begin{align}
\mathbb{E}(\mathfrak{J}) = \mathbb{E}(\mathcal{V}, N, \delta) \bigcap \mathfrak{J},\,\,\,E(\mathfrak{J}) = |\mathbb{E}(\mathfrak{J})|,
\end{align}
where $\mathfrak{J}$ is any subset of $[0, N]$. By the definition of $\mathcal{V}$ and $\mathbb{E}(\mathcal{V}, N, \delta)$, it is easy to see that
\begin{align}
{E}(\mathcal{V}, N, \delta) \ll E([N^{\frac{359}{378}}, N]) + N^{\frac{359}{378}}.
\end{align}
So it remains to estimate $E([N^{\frac{359}{378}}, N])$. Let $N^{\frac{359}{378}} \leqslant X \leqslant N$. Due to the standard dyadic argument, we shall focus on  those $v$ satisfying $\frac{1}{2}X \leqslant v \leqslant X$. Generally, we can consider $2^{-j}X \leqslant v \leqslant 2^{1-j}X \,(j=1,2,...)$, and obtain the desired bound for the exceptional set.

In Sections 3-5, we will estimate the integrals over $\mathfrak{M}$, $\mathfrak{m}$ and $\mathfrak{t}$ and then complete the proof of the Theorem in Section 6.
Now we state the Lemmas required in this paper.\\

{\bf Lemma 2.1}.  {\it  We have}
\begin{align*}
&\mbox{i)}\,\, S_{j}(\alpha) \ll X^{\frac{1}{j}},\,\,\,\,\,T_{j}(\alpha) \ll X^{\frac{1}{j} - 1}\min(X, |\alpha|^{-1});\\
&\mbox{ii)}\,\,S_{2}(\lambda_1 \alpha)  = \kappa L^{-1}T_{2}(\lambda_1 \alpha) + O\big(X^{\frac12}L^{-2}(1+|\alpha|X)\big).
\end{align*}

{\bf Proof}. i) follows from the first derivative estimate for trigonometric integers (see Lemma 4.2 in Titchmarsh [14]). And ii) follows from Abel's summation formula and (2.3).\\

{\bf Lemma 2.2}.  {\it  Suppose that $\alpha$ is a real number, and there exist $a \in \mathbb{Z}$ and $q \in \mathbb{N}$ with
\begin{align*}
(a,q) = 1, \,\,1\leqslant q \leqslant P^{\frac32},\,\,|q\alpha - a|<P^{-\frac32}.
\end{align*}
Then one has}
\begin{align*}
\sum_{p \sim P}(\log p)e(\alpha p^3) \ll P^{\frac{11}{12}+\varepsilon} + \frac{P^{1+\varepsilon}}{q^{\frac12}\big(1+P^3\big|\alpha-\frac{a}{q}\big|\big)^{\frac12}}.
\end{align*}

{\bf Proof}. See Lemma 2.3 in Zhao [19].\\

{\bf Lemma 2.3}. {\it Let $m(2)= \frac{1}{14}$, $m(3)= \frac{1}{36}$ , $\lambda$ and $\mu$ be nonzero constants. For $i=2 \,\mbox{or} \,\,3$, suppose that $X^{\frac{1}{i} - m(i) + \varepsilon} \leqslant Z_i \leqslant X^{\frac {1}{i}}$. Define}
 \begin{align*}
\mathfrak{N}(Z_i) = \{\alpha \in \mathbb{R}:  Z_i \leqslant |S_i(\lambda \alpha)| \leqslant X^{\frac{1}{i}}\}.
\end{align*}
\mbox{i)}\,{\it For $\alpha \in \mathfrak{N}(Z_i)$, there exist coprime integers $a_i$, $q_i$ satisfying}
\begin{align*}
1 \leqslant q_{i} \ll \bigg(\frac{X^{\frac {1}{i} + \varepsilon}}{Z_i}\bigg)^{2},\,\,\,\,\,|\alpha \lambda_{i-1} q_i - a_i| \ll X^{-1}\bigg(\frac{X^{\frac {1}{i}+ \varepsilon}}{Z_i}\bigg)^{2}.
\end{align*}
\mbox{ii)}\,{\it For $j \geqslant 3$, we have}
\begin{align*}
\int_{\mathfrak{N}(Z_i)}|S_i(\lambda \alpha)|^2|S_j(\mu \alpha)|^2K(\alpha)d\alpha \ll \tau(X^{\frac{2}{i} + \frac{2}{j}-1 +\varepsilon} + X^{\frac{4}{i} + \frac{1}{j}-1+\varepsilon}Z_{i}^{-2}).
\end{align*}

{\bf Proof}.
For $i=2$, it follows from Lemma 1 in Wang and Yao [16] that there are coprime integers $a_2$, $q_2$ satisfying
\begin{align}
1 \leqslant q_{2} \ll \bigg(\frac{X^{\frac {1}{2} + \varepsilon}}{Z_2}\bigg)^{4},\,\,\,\,\,|\alpha \lambda_{1} q_2 - a_2| \ll X^{-1}\bigg(\frac{X^{\frac {1}{2}+ \varepsilon}}{Z_2}\bigg)^{4}.
\end{align}
Moreover, we deduce from Lemma 5.6 in Kumchev [10] with $z= X^{\frac{3}{28}}$ and the hypothesis $Z_2 \geqslant X^{\frac12 - m(2) + \varepsilon}$ that
\begin{align}
X^{\frac12 - m(2) + \varepsilon} \leqslant Z_2 \leqslant |S_2(\lambda_1\alpha)|\,&\ll \frac{X^{\frac12 + \varepsilon}}{(q_2+X|\alpha \lambda_{1} q_2 - {a_2}|)^{\frac12}} + X^{\frac17 + \frac{11}{40} + \varepsilon} + X^{\frac{11}{28}+\varepsilon}\notag\\
& \ll \frac{X^{\frac12 + \varepsilon}}{(q_2+X|\alpha \lambda_{1} q_2 - {a_2}|)^{\frac12}}.
\end{align}
Hence, by (2.11), we get
\begin{align*}
1 \leqslant q_{2} \ll \bigg(\frac{X^{\frac {1}{2} + \varepsilon}}{Z_2}\bigg)^{2},\,\,\,\,\,|\alpha \lambda_{1} q_2 - a_2| \ll X^{-1}\bigg(\frac{X^{\frac {1}{2}+ \varepsilon}}{Z_2}\bigg)^{2}.
\end{align*}

For $i=3$, let $P= X^{\frac13}$ and $Q= X^{\frac32}$. By Dirichlet's approximation theorem, there exists two coprime integers $a_3, q_3$ with $1\leqslant q_3 \leqslant Q$ and $|\alpha \lambda_{2} q_3 - a| \leqslant Q^{-1}$. Then it follows from Lemma 2.2 and the hypothesis $Z_3 \geqslant X^{\frac13 - m(3) + \varepsilon}$, we have
\begin{align*}
X^{\frac13 - m(3) + \varepsilon} &\leqslant Z_3 \leqslant |S_3(\lambda_2\alpha)|\ll X^{\frac13 - m(3) + \varepsilon} + \frac{X^{\frac13 + \varepsilon}}{(q_3+X|\alpha \lambda_{2} q_3 - {a_3}|)^{\frac12}}.
\end{align*}
Thus we have
\begin{align}
1 \leqslant q_{3} \ll \bigg(\frac{X^{\frac {1}{3} + \varepsilon}}{Z_3}\bigg)^{2},\,\,\,\,\,|\alpha \lambda_{2} q_3 - a_3| \ll X^{-1}\bigg(\frac{X^{\frac {1}{3}+ \varepsilon}}{Z_3}\bigg)^{2},
\end{align}
which completes the proof of Lemma 2.3 i).

Next we give the proof of Lemma 2.3 ii). By the substitution of variables, it suffices to prove in the case $\lambda =1$. Hence we suppose that $\lambda =1$ in the following proof. Let $Q_i = X^{\frac{2}{i}+\varepsilon}Z_i^{-2}$. We set
\begin{align*}
\mathfrak{N}_i^{'}(q,a) &= \left\{ \alpha: \left|\alpha- \frac{a}{q}\right| \leqslant \frac{Q_i}{qX} \right\}.
\end{align*}
Let $V_i(\alpha)$ be the function of period 1, and defined for $\alpha \in (Q_iX^{-1}, 1+Q_iX^{-1}]$ by
\begin{equation*}
V_i(\alpha) =\left\{\begin{array}{ll}
(q + X|q\alpha - a|)^{-1}, & \alpha \in {\mathfrak{N}_i^{'}}; \\
\\
0, & \alpha \in (Q_iX^{-1}, 1+Q_iX^{-1}]\setminus{\mathfrak{N}_i^{'}},
\end{array}
\right.
\end{equation*}
where $\mathfrak{N}_i^{'}$ is the union of all intervals $\mathfrak{N}_i^{'}(q,a)$ with $1\leqslant a \leqslant q \leqslant Q_i$ and $(a,q)=1$. Let $\mathfrak{N}_i^{*} = \mathfrak{N}_i^{'} + \mathbb{Z}$ be the union of all intervals $\mathfrak{N}_i^{'}(q,a)$ with $1 \leqslant q \leqslant Q_i, a\in \mathbb{Z}$ and $(a,q)=1$. Then we have
\begin{align}
{\mathfrak{N}(Z_i)} \subseteq \mathfrak{N}_i^{*}.
\end{align}
When $\alpha \in \mathfrak{N}(Z_i)$, by combining (2.11)-(2.13), we find that
\begin{align*}
|S_i(\lambda\alpha)| \ll X^{\frac{1}{i}+\varepsilon}V_i(\alpha)^{\frac12}.
\end{align*}
Write
\begin{align*}
|S_j(\mu \alpha)|^2 = \sum_{v}\psi(v)e(\alpha v),
\end{align*}
where
\begin{align*}
\psi(v) = \sum_{{p_1,p_2 \in I_j\atop{\mu(p_1^j - p_2^j) = v}}}(\log p_1)(\log p_2).
\end{align*}
We deduce from [2, Lemma 3] that
\begin{align}
&\int_{\mathfrak{N}(Z_i)}|S_i(\lambda \alpha)|^2|S_j(\mu \alpha)|^2K(\alpha)d\alpha \notag\\
\ll &X^{\frac{2}{i}+\varepsilon}\int_{\mathfrak{N}_i^{*}}V_i(\alpha)\bigg(\sum_{v}\psi(v)e(\alpha v)\bigg)K(\alpha)d\alpha \notag\\
\ll &\tau X^{\frac{2}{i}-1+\varepsilon}(1+\tau)^{1+\varepsilon}\biggr(\sum_{v}\psi(v) + Q_i\sum_{|v|\leqslant \tau}\psi(v)\biggr).
\end{align}
It is easy to see that
\begin{align}
\sum_{v}\psi(v) = |S_j(0)|^2 \ll X^{\frac{2}{j}}.
\end{align}
Note that $\tau \rightarrow 0 $ as $X \rightarrow \infty$. Then it follows from the Prime Number Theorem that
\begin{align}
\sum_{|v|\leqslant \tau}\psi(v) =  \sum_{p\in I_j}(\log p)^2 \ll X^{\frac{1}{j} + \varepsilon}.
\end{align}
Lemma 2.3 ii) can be proved easily by inserting (2.15) and (2.16) into (2.14).\\

{\bf Lemma 2.4}. {\it If $\alpha$ is a real number satisfying that there exists $a \in \mathbb{Z}$ and $q \in \mathbb{N}$ with $(a,q)=1, 1\leqslant q \leqslant P^{\frac34}$ and $|q\alpha - a|\leqslant P^{-\frac94}$, we have
\begin{align*}
\sum_{x \sim P}e(\alpha x^3) \ll \frac{w_3(q)P}{1 + P^3|\alpha - a/q|},
\end{align*}
otherwise we have $\sum\limits_{x \sim P}e(\alpha x^3) \ll P^{\frac34 + \varepsilon}$,
where the multiplicative function $w_3(q)$ is defined by}
\begin{equation}
w_3(p^{3u+v}) =\left\{\begin{array}{ll}
3p^{-u-\frac12}, & u\geqslant 0, v=1; \\
\\
p^{-u-1}, & u\geqslant 0, 2\leqslant v \leqslant 3.
\end{array}
\right.
\end{equation}

{\bf Proof}. See Lemma 2.3 in Zhao [18].\\

{\bf Lemma 2.5}. {\it Let $c$ be a constant and $A$ be a positive constant. For $Q\geqslant 2$, we have}
\begin{align*}
\sum_{1\leqslant q \leqslant Q}d(q)^{c}w_3(q)^2 \ll (\log Q)^{A}.
\end{align*}

{\bf Proof}. See Lemma 2.1 in Zhao [18].\\

{\bf Lemma 2.6}. {\it Suppose that $F \in \{{S_{2}^{4}},\,{S_{4}^{16}},\,{S_{5}^{32}},\,{S_{2}^{2}S_{3}^{4}},\,{S_{2}^{2}S_{4}^{4}},\,{S_{2}^{2}S_{5}^{6}}$\},
then we have}
\begin{align*}
\int_{-1}^{1}|F(\alpha)|d\alpha \ll X^{-1}(F(0))^{1+\varepsilon},\,\,\,\,\int_{-\infty}^{+\infty}|F(\alpha)|K(\alpha)d\alpha \ll \tau X^{-1}(F(0))^{1+\varepsilon}.
\end{align*}

{\bf Proof}. All of these  follow from Vaughan [15] by using Hua's Lemma and the recent result [1] of Vinogradov's mean value theorem. We can also find this result in Lemma 3 in Gao and Liu [5].

\section{The major arc $\mathfrak{M}$}
\setcounter{equation}{0}
We subdivide the major arc $\mathfrak{M}$ into two parts: $\mathfrak{M} = {\mathfrak{M}_{1}} \bigcup {\mathfrak{M}_{2}}$, where
\begin{align*}
\mathfrak{M}_1 = \{\alpha: |\alpha| \leqslant X^{-\frac{11}{12}}\}.
\end{align*}
Thus we have
\begin{align}
I(\tau, v, \mathfrak{M})= I(\tau, v, \mathfrak{M}_1) + I(\tau, v, \mathfrak{M}_2).
\end{align}
Let $H(\alpha) = \kappa L^{-1}\prod\limits_{j=1}^{4}T_{j+1}(\lambda_j\alpha)e(-\alpha v)K(\alpha)$.
It follows that
\begin{align}
&I(\tau, v, \mathfrak{M}_1) - \int_{\mathfrak{M}_1}H(\alpha)d\alpha\notag\\
= &\int_{\mathfrak{M}_1}\left(S_2(\lambda_1\alpha)-\kappa L^{-1}T_2(\lambda_1\alpha)\right)\prod_{j=2}^{4}T_{j+1}(\lambda_j\alpha)e(-\alpha v)K(\alpha)d\alpha\notag\\
&\,\,+ \int_{\mathfrak{M}_1}S_2(\lambda_1\alpha)\left(S_3(\lambda_2\alpha)- T_3(\lambda_2\alpha)\right)\prod_{j=3}^{4}T_{j+1}(\lambda_j\alpha)e(-\alpha v)K(\alpha)d\alpha\notag\\
&\,\,+ \int_{\mathfrak{M}_1}\prod_{j=1}^{2}S_{j+1}(\lambda_j\alpha)\left(S_4(\lambda_3\alpha)- T_4(\lambda_3\alpha)\right)T_{5}(\lambda_4\alpha)e(-\alpha v)K(\alpha)d\alpha\notag\\
&\,\,+ \int_{\mathfrak{M}_1}\prod_{j=1}^{3}S_{j+1}(\lambda_j\alpha)\left(S_5(\lambda_4\alpha)- T_5(\lambda_4\alpha)\right)e(-\alpha v)K(\alpha)d\alpha\notag\\
=: &J_1 +J_2+J_3+J_4.
\end{align}
Since the computations for $J_2, J_3$ and $J_4$ are similar to the corresponding ones in Section 3 in Mu and Gao [12]. We have
\begin{align}
J_i \ll \tau^2 X^{\frac{17}{60}}L^{-1}.\,\,\,\,\,(i=2,3,4)
\end{align}
Next, we restrict our attention to estimating $J_1$. Note that
\begin{align}
J_1 \ll \tau^2\int_{\mathfrak{M}_1}\left|S_2(\lambda_1\alpha)-\kappa L^{-1}T_2(\lambda_1\alpha)\right|\prod_{j=2}^{4}|T_{j+1}(\lambda_j\alpha)|d\alpha.
\end{align}
It follows from partial summation  that
\begin{align}
S_2(\lambda_1\alpha) = \int_{I_2}e(\lambda_1t^2\alpha)d\bigg(\sum_{m\leqslant t, m\in I_2}\rho(m)\bigg).
\end{align}
Then by Lemma 2.1 and (3.4), we obtain
\begin{align}
J_1 &\ll \tau^2X^{\frac12}L^{-2}\int_{0}^{\frac{1}{X}}\prod_{j=2}^{4}|T_{j+1}(\lambda_j\alpha)|d\alpha + \tau^2X^{\frac32}L^{-2}\int_{\frac{1}{X}}^{\frac{P}{X}}\alpha\prod_{j=2}^{4}|T_{j+1}(\lambda_j\alpha)|d\alpha\notag\\
&\ll \tau^2 X^{\frac{17}{60}}L^{-1}.
\end{align}
By the argument similar to Section 3 in Mu and Gao [12], we can obtain
\begin{align}
\int_{\mathfrak{M}_1}H(\alpha)d\alpha \gg  \tau^2 X^{\frac{17}{60}}L^{-1},
\end{align}
and
\begin{align}
|I(\tau, v, \mathfrak{M}_2)| \ll \tau^2 X^{\frac{17}{60}}L^{-1}.
\end{align}
This in combination with (3.2), (3.3), (3.6)-(3.8) gives that
\begin{align}
|I(\tau, v, \mathfrak{M})| \gg \tau^2 X^{\frac{17}{60}}L^{-1}.
\end{align}

\section{The Trivial arc $\mathfrak{t}$}
\setcounter{equation}{0}
It follows from Section 4 in Mu and Gao [12] that
\begin{align}
|I(\tau, v,\mathfrak{t})| \ll \tau^2X^{\frac{17}{60} - \varepsilon}.
\end{align}

\section {The Minor Arc ${\mathfrak{m}}$}
\setcounter{equation}{0}
For convenience, we write
\begin{align}
\sigma = \frac{1}{378}.
\end{align}
We divide ${\mathfrak{m}}$ into four disjoint parts:
\begin{align}
{\mathfrak{m}} = {\mathfrak{m}_1} \bigcup {\mathfrak{m}_2} \bigcup {\mathfrak{m}_3} \bigcup {\mathfrak{m}_4},
\end{align}
where
\begin{align*}
{\mathfrak{m}_1} &= \{\alpha \in \mathfrak{m}: |S_2(\lambda_1\alpha)|< X^{\frac12 - 27\sigma + 2\varepsilon},\,|S_3(\lambda_2\alpha)|< X^{\frac13 - \frac{21}{2}\sigma + 2\varepsilon}\},\\
{\mathfrak{m}_2} &= \{\alpha \in \mathfrak{m}: |S_2(\lambda_1\alpha)|< X^{\frac12 - 27\sigma + 2\varepsilon},\,|S_3(\lambda_2\alpha)|\geqslant X^{\frac13 - \frac{21}{2}\sigma + 2\varepsilon}\},\\
{\mathfrak{m}_3} &= \{\alpha \in \mathfrak{m}: |S_2(\lambda_1\alpha)|\geqslant X^{\frac12 - 27\sigma + 2\varepsilon},\,|S_3(\lambda_2\alpha)|< X^{\frac13 - \frac{21}{2}\sigma  + 2\varepsilon}\},\\
{\mathfrak{m}_4} &= \{\alpha \in \mathfrak{m}: |S_2(\lambda_1\alpha)|\geqslant X^{\frac12 - 27\sigma + 2\varepsilon},\,|S_3(\lambda_2\alpha)|\geqslant X^{\frac13 - \frac{21}{2}\sigma  + 2\varepsilon}\}.
\end{align*}

\subsection {The estimation over ${\mathfrak{m}_1}$}
\hspace{4.6mm} {\bf Lemma 5.1}. {\it For $\alpha \in \mathfrak{m}_1$, we have}
\begin{align*}
\int_{\mathfrak{m}_1}|S_2(\lambda_1 \alpha)|^2|S_3(\lambda_2 \alpha)|^{12}K(\alpha)d\alpha \ll \tau X^{4- 114\sigma + \varepsilon}.
\end{align*}

{\bf Proof}. The argument is analogous to Lemma 7.1 in Ge and Zhao [6]. For convenience, we write $G_t(\alpha) = |S_2(\lambda_1 \alpha)|^2|S_3(\lambda_2 \alpha)|^{t-2}$ and $J_t = \int_{\mathfrak{m}_1}|S_2(\lambda_1 \alpha)|^2|S_3(\lambda_2 \alpha)|^{t}K(\alpha)d\alpha$. Therefore,
\begin{align*}
J_t &= \int_{\mathfrak{m}_1}S_3(\lambda_2 \alpha)S_3(-\lambda_2 \alpha)G_t(\alpha)K(\alpha)d\alpha\\
&= \sum_{p \in I_3}(\log p)\int_{\mathfrak{m}_1}e(\alpha \lambda_2p^3)S_3(-\lambda_2 \alpha)G_t(\alpha)K(\alpha)d\alpha \\
&\ll L\sum_{n \in I_3}\left|\int_{\mathfrak{m}_1}e(\alpha \lambda_2n^3)S_3(-\lambda_2 \alpha)G_t(\alpha)K(\alpha)d\alpha\right|.
\end{align*}
An application of Cauchy's inequality reveals that
\begin{align}
J_t \ll LX^{\frac16}{\mathfrak{J}_t}^{\frac12},
\end{align}
where
\begin{align}
{\mathfrak{J}_t} = \sum_{n \in I_3}\left|\int_{\mathfrak{m}_1}e(\alpha \lambda_2n^3)S_3(-\lambda_2 \alpha)G_t(\alpha)K(\alpha)d\alpha\right|^2.
\end{align}
For the sum ${\mathfrak{J}_t}$, we have
\begin{align}
{\mathfrak{J}_t} &= \sum_{n \in I_3}\int_{\mathfrak{m}_1}\int_{\mathfrak{m}_1}|G_t(\alpha)G_t(-\beta)|S_3(-\lambda_2 \alpha)S_3(\lambda_2 \beta)e(\lambda_2n^3(\alpha - \beta))K(\alpha)K(\beta)d\alpha d\beta\notag\\
&\leqslant \int_{\mathfrak{m}_1}|G_t(\beta)S_3(\lambda_2 \beta)|F(\beta)K(\beta)d\beta,
\end{align}
where
\begin{align*}
T(x) = \sum_{n \in I_3}e(xn^3),
\end{align*}
and
\begin{align}
F(\beta) = \int_{\mathfrak{m}_1}|G_t(\alpha)S_3(-\lambda_2 \alpha)T(\lambda_2(\alpha-\beta))|K(\alpha)d\alpha.
\end{align}
Let $\mathcal{M}_\beta(r,b) = \{\alpha \in {\mathfrak{m}_1} : |r\lambda_2(\alpha-\beta) -b| \leqslant X^{-\frac34}\}$. Then the set $\mathcal{M}_\beta(r,b) \neq \emptyset$ forces that
\begin{align*}
|b + r\lambda_2\beta| \leqslant |r\lambda_2(\alpha-\beta) -b| +  |r\lambda_2\alpha|   \leqslant X^{-\frac34} + r|\lambda_2|R.
\end{align*}
Let $\mathcal{B} = \{b \in \mathbb{Z}:  |b + r\lambda_2\beta|\leqslant X^{-\frac34} + r|\lambda_2|R\}$. We divide the set $\mathcal{B}$ into two disjoint sets $\mathcal{B}_1 = \{b \in \mathbb{Z}:  |b + r\lambda_2\beta|\leqslant r|\lambda_2|\tau^{-1}\}$ and $\mathcal{B}_2 = \mathcal{B}\setminus \mathcal{B}_1.$ Let
\begin{align*}
\mathcal{M}_\beta = \bigcup_{1\leqslant r\leqslant X^{\frac14}}\bigcup_{{b \in \mathcal{B}}\atop{(b,r)=1}}\mathcal{M}_\beta(r,b).
\end{align*}
Then by Lemma 2.4 and Cauchy's inequality, we have
\begin{align}
F(\beta) &\ll X^{\frac13}\int_{{\mathcal{M}_\beta}\bigcap {\mathfrak{m}_1}}\frac{|G_t(\alpha)S_3(-\lambda_2 \alpha)|w_3(r)K(\alpha)}{1 + X|\lambda_2(\alpha - \beta) - b/r|}d\alpha + X^{\frac 14 + \varepsilon}\int_{\mathfrak{m}_1}|G_t(\alpha)S_3(-\lambda_2 \alpha)|K(\alpha)d\alpha\notag\\
& \ll X^{\frac13}\left(\int_{{\mathfrak{m}_1}}|G(\alpha)|^2 K(\alpha)d\alpha\right)^{\frac12}L(\beta)^{\frac12}  + X^{\frac 14 + \varepsilon}\int_{\mathfrak{m}_1}|G_t(\alpha)S_3(-\lambda_2 \alpha)|K(\alpha)d\alpha,
\end{align}
where
\begin{align}
L(\beta) = \int_{{\mathcal{M}_\beta}}\frac{|S_3(-\lambda_2 \alpha)|^2w_3(r)^2K(\alpha)}{(1 + X|\lambda_2(\alpha - \beta) - b/r|)^2}d\alpha.
\end{align}
Furthermore, we divide $L(\beta)$ into two parts:
\begin{align}
L(\beta) &= \sum_{j=1,2}\sum_{1\leqslant r\leqslant X^{\frac14}}\sum_{{b \in \mathcal{B}_j}\atop{(b,r)=1}}\int_{{\mathcal{M}_\beta(r,b)}}\frac{|S_3(-\lambda_2 \alpha)|^2w_3(r)^2K(\alpha)}{(1 + X|\lambda_2(\alpha - \beta) - b/r|)^2}d\alpha \notag\\
&=: L_1(\beta) + L_2(\beta).
\end{align}
For $L_1(\beta)$, we have
\begin{align*}
L_1(\beta) &\ll \tau^2\sum_{1\leqslant r\leqslant X^{\frac14}}w_3(r)^2\sum_{{b \in \mathcal{B}_1}\atop{(b,r)=1}}\int_{|r\lambda_2\gamma| \leqslant X^{-\frac34}}\frac{|S_3(\lambda_2(\beta+\gamma))+b/r|^2}{(1 + X|\lambda_2\gamma|)^2}d\gamma \\
&\ll \tau^2\sum_{1\leqslant r\leqslant X^{\frac14}}w_3(r)^2\int_{|r\lambda_2\gamma| \leqslant X^{-\frac34}}\frac{U({\mathcal{B}}^{*}_1)}{(1 + X|\lambda_2\gamma|)^2}d\gamma,
\end{align*}
where
\begin{align*}
U({\mathcal{B}}^{*}_1) &= \sum_{b \in {\mathcal{B}}^{*}_1}|S_3(\lambda_2(\beta+\gamma))+b/r|^2,\\
{\mathcal{B}}^{*}_1 &= \{b \in \mathbb{Z}: -r([|\lambda_2|\tau^{-1}]+1) < b + r\lambda_2\beta \leqslant r([|\lambda_2|\tau^{-1}]+1)\}.
\end{align*}
Then we get
\begin{align*}
U({\mathcal{B}}^{*}_1) &= \sum_{p_1,p_2 \in I_3}\sum_{b \in {\mathcal{B}}^{*}_1}e((\lambda_2(\beta+\gamma) + b/r)(p_1^3 - p_2^3))\\
&\leqslant \sum_{p_1,p_2 \in I_3}\biggr|\sum_{b \in {\mathcal{B}}^{*}_1}e(b/r(p_1^3 - p_2^3))\biggr|\\
&= 2r(|\lambda_2|\tau^{-1} +1)\sum_{{p_1,p_2 \in I_3}\atop{p_1^3-p_2^3 \equiv 0 (\bmod r)}}1\\
&\ll r\tau^{-1}X^{\frac23}r^{-2}\sum_{{1\leqslant b_1,b_2<r, (b_1b_2, r)=1}\atop{b_1^3-b_2^3 \equiv 0 (\bmod r)}}1\\
&\ll \tau^{-1}X^{\frac23}\sum_{{1\leqslant b<r}\atop{b^3 \equiv 1 (\bmod r)}}1    \,\,\,\,\,\ll \tau^{-1}X^{\frac23}d(r)^{c}.
\end{align*}
Thus we deduce from Lemma 2.5 that
\begin{align}
L_1(\beta) &\ll \tau X^{\frac 23}\sum_{1 \leqslant r \leqslant X^{\frac14}}w_3(r)^2d(r)^{c}\int_{|r\lambda_2\gamma| < X^{-\frac34}}\frac{1}{(1+X|\lambda_2\gamma|)^2}d\gamma\notag\\
&\ll \tau X^{-\frac13}\sum_{1 \leqslant r \leqslant X^{\frac14}}w_3(r)^2d(r)^{c}\,\,\ll \tau X^{-\frac13 + \varepsilon}.
\end{align}
Next, we give the estimation of $L_2(\beta)$. Without loss of generality, it suffices to consider the set
\begin{align*}
{\mathcal{B}}^{'}_2 = \{b \in \mathbb{Z}: r|\lambda_2|\tau^{-1} < b + r\lambda_2\beta \leqslant  r|\lambda_2|R + X^{-\frac34}\}
\end{align*}
which falls in the set
\begin{align*}
{\mathcal{B}}^{*}_2 = \{b \in \mathbb{Z}: r\kappa_1 < b + r\lambda_2\beta \leqslant  r\kappa_2\},
\end{align*}
where $\kappa_1 = [|\lambda_2|\tau^{-1}]$ and $\kappa_2 = [|\lambda_2|R] + 2$. Then we get
\begin{align*}
L_2(\beta) &\ll \sum_{1\leqslant r\leqslant X^{\frac14}}\sum_{b \in \mathcal{B}^{*}_2}\int_{{\mathcal{M}_\beta(r,b)}}\frac{|S_3(-\lambda_2 \alpha)|^2w_3(r)^2K(\alpha)}{(1 + X|\lambda_2(\alpha - \beta) - b/r|)^2}d\alpha\\
&\ll \sum_{1\leqslant r\leqslant X^{\frac14}}w_3(r)^2\sum_{b \in \mathcal{B}^{*}_2}\int_{{\mathcal{M}_\beta(r,b)}}\frac{|S_3(-\lambda_2 \alpha)|^2|\alpha|^{-2}}{(1 + X|\lambda_2(\alpha - \beta) - b/r|)^2}d\alpha\\
&\ll \sum_{1\leqslant r\leqslant X^{\frac14}}w_3(r)^2\sum_{\kappa_1 \leqslant k<\kappa_2}\frac{1}{(k-1)^2}\sum_{rk < b+r\lambda_2\beta \leqslant r(k+1)}\int_{{\mathcal{M}_\beta(r,b)}}\frac{|S_3(-\lambda_2 \alpha)|^2}{(1 + X|\lambda_2(\alpha - \beta) - b/r|)^2}d\alpha\\
&\ll \sum_{1\leqslant r\leqslant X^{\frac14}}w_3(r)^2\sum_{\kappa_1 \leqslant k<\kappa_2}\frac{1}{(k-1)^2}\int_{|r\lambda_2\gamma|\leqslant X^{-\frac34}}\frac{U(\mathcal{C}_k)}{(1+X|\lambda_2\gamma|)^2}d \gamma,
\end{align*}
where $\mathcal{C}_k = \{b \in \mathbb{Z}: rk < b+r\lambda_2\beta \leqslant r(k+1)\}.$ Similar to the estimation of $U(\mathcal{B}^{*}_1)$, we have $U(\mathcal{C}_k) \ll X^{\frac23}d(r)^c$. Thus we have
\begin{align}
L_2(\beta) \ll X^{-\frac13}\sum_{1\leqslant r\leqslant X^{\frac14}}w_3(r)^2d(r)^c\sum_{\kappa_1 \leqslant k<\kappa_2}\frac{1}{(k-1)^2} \,\,\ll \tau X^{-\frac13+\varepsilon}.
\end{align}
On combining (5.7)-(5.11), we have
\begin{align}
F(\beta) \ll X^{\frac16+\varepsilon}\left(\int_{{\mathfrak{m}_1}}|G_t(\alpha)|^2 K(\alpha)d\alpha\right)^{\frac12}\tau^{\frac12}  + X^{\frac 14 + \varepsilon}\int_{\mathfrak{m}_1}|G_t(\alpha)S_3(-\lambda_2 \alpha)|K(\alpha)d\alpha.
\end{align}
Consequently, by (5.3),(5.5) and (5.12), we obtain
\begin{align}
J_t \ll& \tau^{\frac14}X^{\frac14+\varepsilon}\left(\int_{{\mathfrak{m}_1}}|G_t(\alpha)|^2 K(\alpha)d\alpha\right)^{\frac14}{\left(\int_{\mathfrak{m}_1}|G_t(\alpha)S_3(-\lambda_2 \alpha)|K(\alpha)d\alpha\right)}^{\frac12}\notag\\
&+ X^{\frac{7}{24}+\varepsilon}\int_{\mathfrak{m}_1}|G_t(\alpha)S_3(-\lambda_2 \alpha)|K(\alpha)d\alpha.
\end{align}
From the definition of ${\mathfrak{m}_1}$, we find that
\begin{align}
\int_{{\mathfrak{m}_1}}|G_t(\alpha)|^2 K(\alpha)d\alpha &\ll \sup\limits_{\alpha \in  {\mathfrak{m}_1}}|S_2(\lambda_1\alpha)^2S_3(\lambda_2\alpha)^{t-4}|J_{t}\notag\\
&\ll X^{\frac13 t -\frac13 - \left(12\sigma + \frac{21}{2}\sigma t\right)+\varepsilon}J_t.
\end{align}
Moreover, an application of  Cauchy's inequality yields the estimate
\begin{align}
\int_{\mathfrak{m}_1}|G_t(\alpha)S_3(-\lambda_2 \alpha)|K(\alpha)d\alpha \ll {J_t}^{\frac12}{J_{t-2}}^{\frac12}.
\end{align}
Concluding (5.13)-(5.15), we have
\begin{align}
J_{t} \ll \tau^{\frac12}{J^{\frac12}_{t-2}}X^{\frac13 +\frac16t - \left(6\sigma + \frac{21}{4}\sigma t\right) + \varepsilon} + X^{\frac{7}{12} + \varepsilon}J_{t-2}.
\end{align}
It follows from Lemma 2.6 that
\begin{align}
J_4 \ll \int_{-\infty}^{\infty}|S_2(\lambda_1\alpha)^2S_3(\lambda_2\alpha)^4|K(\alpha)d\alpha \ll \tau X^{\frac43 + \varepsilon}.
\end{align}
By simple calculation, we can deduce from (5.16)-(5.17) that
\begin{align}
J_6 \ll \tau X^{2-\frac{63}{2}\sigma+\varepsilon},\,\,\,\,J_8 \ll \tau X^{\frac83 - 63\sigma+\varepsilon},\,\,\,\,J_{10} \ll \tau X^{\frac{10}{3} - 90\sigma +\varepsilon},\,\,\,\,J_{12} \ll \tau X^{4-114\sigma+\varepsilon},
\end{align}
which completes the proof of Lemma 5.1.\\

{\bf Proposition 1}. {\it We have}
\begin{align}
\int_{\mathfrak{m}_1}\prod_{j=1}^{4}\big|S_{j+1}(\lambda_j\alpha)\big|^2K(\alpha)d\alpha \ll \tau X^{\frac{47}{30} - 19\sigma + \varepsilon}.
\end{align}

{\bf Proof}. Lemmas 2.6 and 5.1 and H\"{o}lder's inequality imply that
\begin{align*}
&\int_{\mathfrak{m}_1}|S_2(\lambda_1 \alpha)S_3(\lambda_2 \alpha)S_4(\lambda_3 \alpha)S_5(\lambda_4 \alpha)|^2K(\alpha)d\alpha\\
 \ll \,&\left(\int_{\mathfrak{m}_1}|S_2(\lambda_1 \alpha)|^2|S_3(\lambda_2 \alpha)|^{12}K(\alpha)d\alpha\right)^{\frac16}
\left(\int_{-\infty}^{\infty}|S_2(\lambda_1 \alpha)|^2|S_4(\lambda_3 \alpha)|^{4}K(\alpha)d\alpha\right)^{\frac12}\\
&\times \left(\int_{-\infty}^{\infty}|S_2(\lambda_1 \alpha)|^2|S_5(\lambda_4 \alpha)|^{6}K(\alpha)d\alpha\right)^{\frac13}\\
\ll \,& \,\,\tau X^{\frac{47}{30} - 19\sigma + \varepsilon}.
\end{align*}

\subsection {The estimation over ${\mathfrak{m}_2}$}
{\bf Proposition 2}. {\it We have}
\begin{align}
\int_{\mathfrak{m}_2}\prod_{j=1}^{4}\big|S_{j+1}(\lambda_j\alpha)\big|^2K(\alpha)d\alpha\ll \tau X^{\frac{47}{30} - 54\sigma + \varepsilon}.
\end{align}

{\bf Proof}. On recalling the definition of ${\mathfrak{m}_2}$ and applying H\"{o}lder's inequality as well as Lemma 2.3 ii) with $i=3$ and $j=4$, we arrive at the conclusion that
\begin{align*}
&\int_{\mathfrak{m}_1}|S_2(\lambda_1 \alpha)S_3(\lambda_2 \alpha)S_4(\lambda_3 \alpha)S_5(\lambda_4 \alpha)|^2K(\alpha)d\alpha\\
\ll \,& \sup_{\alpha \in {\mathfrak{m}_2}}\big(|S_2(\lambda_1 \alpha)|^2|S_5(\lambda_4 \alpha)|^2\big)\times \left(\int_{\mathfrak{m}_2}|S_3(\lambda_2 \alpha)|^2|S_4(\lambda_3 \alpha)|^{2}K(\alpha)d\alpha\right)\\
\ll \,&  \tau X^{\frac{47}{30} - 54\sigma + \varepsilon}.
\end{align*}

\subsection {The estimation over ${\mathfrak{m}_3}$}
{\bf Proposition 3}. {\it We have}
\begin{align}
\int_{\mathfrak{m}_3}\prod_{j=1}^{4}\big|S_{j+1}(\lambda_j\alpha)\big|^2K(\alpha)d\alpha \ll\tau X^{\frac{47}{30} - 21\sigma + \varepsilon}.
\end{align}

{\bf Proof}. It follows from the definition of ${\mathfrak{m}_3}$, H\"{o}lder's inequality and Lemma 2.3 ii) with $i=2$ and $j=5$ that
\begin{align*}
&\int_{\mathfrak{m}_3}|S_2(\lambda_1 \alpha)S_3(\lambda_2 \alpha)S_4(\lambda_3 \alpha)S_5(\lambda_4 \alpha)|^2K(\alpha)d\alpha\\
\ll \,&\sup_{\alpha \in {\mathfrak{m}_2}}\big(|S_3(\lambda_2 \alpha)|\big)^2\left(\int_{\mathfrak{m}_3}|S_2(\lambda_1 \alpha)|^2|S_5(\lambda_4 \alpha)|^{2}K(\alpha)d\alpha\right)^{\frac57}\left(\int_{-\infty}^{\infty}|S_2(\lambda_1 \alpha)|^4K(\alpha)d\alpha\right)^{\frac17}\\
&\times \left(\int_{-\infty}^{\infty}|S_4(\lambda_3 \alpha)|^{16}K(\alpha)d\alpha\right)^{\frac18}
\left(\int_{-\infty}^{\infty}|S_5(\lambda_4 \alpha)|^{32}K(\alpha)d\alpha\right)^{\frac{1}{56}}\\
\ll \,&\tau X^{\frac{47}{30} - 21\sigma + \varepsilon}.
\end{align*}

\subsection {The estimation over ${\mathfrak{m}_4}$}
{\bf Proposition 4}. {\it We have}
\begin{align}
\int_{\mathfrak{m}_4}\prod_{j=1}^{4}\big|S_{j+1}(\lambda_j\alpha)\big|^2K(\alpha)d\alpha \ll \tau X^{\frac{47}{30} - \frac{359}{378}\left(1-75\sigma\right)+\varepsilon}.
\end{align}

{\bf Proof}.
We divide ${\mathfrak{m}_4}$ into disjoint sets such that for $\alpha \in S(Z_2,Z_3,y)$, define
\begin{align}
\mathscr{A} &= S(Z_2,Z_3,y)\notag\\
&= \{|S_2(\lambda_1 \alpha)| \sim Z_2,\,\,|S_3(\lambda_2 \alpha)| \sim Z_3,\,\,|\alpha| \sim y\},
\end{align}
where $Z_2 = 2^{k_1}X^{\frac12 -27\sigma + \varepsilon}$, $Z_3 = 2^{k_2}X^{\frac13-\frac{21}{2}\sigma + \varepsilon}$, $y = 2^{k_3}X^{-\frac{5}{6}}$ for some non-negative integers $k_1,k_2,k_3$. To prove the first part of Theorem 1.2, let $N$ run through the sequences
\begin{align}
N = q^{\frac{378}{359(1-75\sigma)}}.
\end{align}

Moreover, denote by $m(\mathscr{A})$ the Lebesgue measure of $\mathscr{A}$. On replacing Lemma 5.4 in Mu and Gao [12]  with Lemma 2.2 i) in this paper, we can deduce from the arguments similar to Lemma 5.6 in Mu and Gao [12] that
\begin{align}
m(\mathscr{A}) \ll yX^{\frac23 - \frac{359}{378}\left(1-75\sigma\right)+\varepsilon}(Z_2Z_3)^{-2}.
\end{align}
Then it follows from (5.25) and (2.1) that
\begin{align*}
\int_{\mathscr{A}}|S_2(\lambda_1\alpha)S_3(\lambda_2\alpha)|^2K(\alpha)d\alpha &\ll (Z_2Z_3)^2\cdot\min((\tau^2, y^{-2})\cdot m(\mathscr{A})\\
&\ll \tau X^{\frac23 - \frac{359}{378}\left(1-75\sigma\right)+\varepsilon},
\end{align*}
this together with Lemma 2.1 i) gives
\begin{align}
&\int_{\mathfrak{m}_4}|S_2(\lambda_1 \alpha)S_3(\lambda_2 \alpha)S_4(\lambda_3 \alpha)S_5(\lambda_4 \alpha)|^2K(\alpha)d\alpha\notag\\
\ll \,& L^3X^{\frac{9}{10}}\cdot \max_{\mathscr{A}}\int_{\mathscr{A}}|S_2(\lambda_1\alpha)S_3(\lambda_2\alpha)|^2K(\alpha)d\alpha \notag\\
\ll \,&\tau X^{\frac{47}{30} - \frac{359}{378}\left(1-75\sigma\right)+\varepsilon}.
\end{align}
Therefore, Propositions 1-4 together yield
\begin{align}
\int_{\mathfrak{m}}\prod_{j=1}^{4}\big|S_{j+1}(\lambda_j\alpha)\big|^2K(\alpha)d\alpha \ll \tau X^{\frac{47}{30} - 19\sigma + \varepsilon}.
\end{align}

\section{Proof of the Theorems}
\setcounter{equation}{0}
On recalling (3.9), (4.1) and (5.27), we arrive at the conclusion that
\begin{align}
\left|\sum_{v \in \mathbb{E}([\frac12 X, X])}I(\tau, v, \mathfrak{m})\right| \gg \tau^2X^{\frac {17}{60}}E\bigg(\bigg[\frac12 X, X\bigg]\bigg).
\end{align}
It follows from Cauchy's inequality and (5.27) that
\begin{align}
&\left|\sum_{v \in \mathbb{E}([\frac12 X, X])}I(\tau, v, \mathfrak{m})\right|\notag\\
\ll& \left(\int_{-\infty}^{\infty}\bigg|\sum_{v \in \mathbb{E}([\frac12 X, X])}e(-\alpha v)\bigg|^2K(\alpha)d \alpha \right)^{\frac 12}\left(\int_{\mathfrak{m}}\prod_{j=1}^{4}\big|S_{j+1}(\lambda_j\alpha)\big|^2K(\alpha)d\alpha\right)^{\frac 12}\\
\ll& (\tau X^{\frac{47}{30} - 19\sigma + \varepsilon})^\frac 12\left(\sum_{v_1,v_2 \in \in \mathbb{E}([\frac12 X, X])}\max(0, \tau - |v_1 - v_2|)\right)^\frac 12\notag\\
\ll& \tau{E\bigg(\bigg[\frac12 X, X\bigg]\bigg)}^{\frac 12}(X^{\frac{47}{30} - 19\sigma  + \varepsilon})^\frac 12.
\end{align}
On combining (6.1) and (6.3), we have
\begin{align}
E\bigg(\bigg[\frac12 X, X\bigg]\bigg) \ll \tau^{-2}X^{1- 19\sigma +\varepsilon}.
\end{align}
Recalling that $\tau = X^{-\delta}$ and we denote a well-spaced sequence by $\mathcal{V}$ and $v_i \in \mathcal{V}$ for $i = 1, 2, ...$. When $N$ is defined by (5.24), we can deduce from (6.4) that
\begin{align}
E(\mathcal{V}, N, \delta) &\ll E([N^{\frac{359}{378}}, N]) + N^{\frac{359}{378}}\notag\\
&= \sum_{j=1}^{\big[\frac{19}{378}\log_2 N\big] + 1}E\bigg(\bigg[\frac12 N, N\bigg]\bigg) + N^{\frac{359}{378}}\notag\\
& \ll \sum_{j=1}^{\big[\frac{19}{378}\log_2 N\big] + 1}\left(\frac{N}{2^{j-1}}\right)^{1-\frac{19}{378}+2\delta+\varepsilon} + N^{\frac{359}{378}}\notag\\
& \ll N^{1-\frac{19}{378}+2\delta+\varepsilon}.
\end{align}

Obviously, there are infinitely many $q$ we would have taken since $\lambda_1/\lambda_2$ is irrational, and this gives the sequence $N_j \rightarrow \infty$. This completes the proof of the first part of Theorem 1.2.

Next, we begin to prove the second part of Theorem 1.2. Let $N$ be sufficiently large, $N^{\frac{359}{378}} \leqslant X \leqslant N$. Assume that the convergent denominators $q_j$ for $\lambda_1/\lambda_2$ satisfy (1.5), then we can modify our work in Section 5. Suppose that
\begin{align}
\sigma = \chi,
\end{align}
where $\chi$ is given by (1.7).
An argument similar to that in the proofs of Section 5 shows that for $j=1,2,3$, we obtain
\begin{align}
\int_{\mathfrak{m}_j}\prod_{j=1}^{4}\big|S_{j+1}(\lambda_j\alpha)\big|^2K(\alpha)d\alpha \ll \tau X^{\frac{47}{30} - 19\chi + \varepsilon}.
\end{align}
We know from (1.5) that there is a convergent $a/q$ to $\lambda_1/\lambda_2$ with
\begin{align*}
X^{(1-75\chi)(1-\omega)} \ll q \ll X^{(1-75\chi)}.
\end{align*}
If not, there would be two convergent denominators $q_j$ and $q_{j+1}$ satisfying
\begin{align*}
q_j \ll X^{(1-75\chi)(1-\omega)},\,\,\,\,q_{j+1} \gg X^{(1-75\chi)}.
\end{align*}
Thus we have $q_j \ll q^{1-\omega}_{j+1}$, which is in contradiction with (1.5).
Proceeding as in the proof of Lemma 5.6 in Mu and Gao [12], we have
\begin{align*}
m(\mathscr{A}) \ll yX^{\frac23 - (1-75\chi)(1-\omega)+\varepsilon}(Z_2Z_3)^{-2}.
\end{align*}
Thus the expression corresponding to (5.26) is
\begin{align}
\int_{\mathfrak{m}_4}\prod_{j=1}^{4}\big|S_{j+1}(\lambda_j\alpha)\big|^2K(\alpha)d\alpha &\ll \tau X^{\frac{47}{30} - (1-75\chi)(1-\omega) + \varepsilon}\notag\\
&\ll \tau X^{\frac{47}{30} - 19\chi + \varepsilon}
\end{align}
by our choice of $\chi$. Then the expression corresponding to (5.27) is
\begin{align}
\int_{\mathfrak{m}}\prod_{j=1}^{4}\big|S_{j+1}(\lambda_j\alpha)\big|^2K(\alpha)d\alpha \ll \tau X^{\frac{47}{30} - 19\chi + \varepsilon}.
\end{align}
Furthermore, on taking account of (6.1) and (6.2), we conclude that
\begin{align}
E\bigg(\bigg[\frac12 X, X\bigg]\bigg) \ll \tau^{-2}X^{1 - 19\chi + \varepsilon} \ll X^{1 - 19\chi + 2\delta + \varepsilon}.
\end{align}
By the same argument as in (6.5), the proof of the second part of Theorem 1.2 is complete.

{\bf{Acknowledgement}}. The author would like to thank the anonymous referee for his/her patience and time in refereeing this manuscript.

\end{document}